\documentclass[12pt]{article}
     \usepackage{amsmath,amsthm} 
     \usepackage{amssymb} 
     \usepackage{enumerate} 
\usepackage[all]{xy}

\newcommand{\OPN}{\mathcal O_{\mathbb P^N}}
\newcommand{\FF}{\mathcal F}
\newcommand{\stir}[2]{\genfrac{[}{]}{0pt}{}{#1}{#2}}

\newtheorem{Theorem}{\quad Theorem}[section]

\newtheorem{Corollary}[Theorem]{\quad Corollary} 

\newtheorem{Lemma}[Theorem]{\quad Lemma}

\newtheorem{remark}[Theorem]{\quad Remark}

    \newenvironment{pf}{ \noindent {\bf Proof:} }{\hfill $\diamond $ \vspace{0.5cm}}

     \title{The Euler characteristic as a polynomial in the Chern classes} 
\author{ Cristina Bertone\\ 
           Dipartimento di Matematica dell'Universit\`{a} \\ 
           Via Carlo Alberto 10 \\ 
           10123 Torino, Italy \\ {\small cristina.bertone@unito.it}            }

\date{}

\begin{document}
\maketitle


\begin{abstract}
In this paper we obtain some explicit expressions for the Euler characteristic of a rank $n$ coherent  sheaf $\FF$ on $\mathbb P^N$ and of its twists $\FF(t)$ as polynomials in the Chern classes $c_i(\FF)$, also giving algorithms for the computation.
The employed methods  use techniques of umbral calculus involving symmetric functions and Stirling numbers.
\end{abstract}

\section*{Introduction}

The aim of this paper is to find a general polynomial expression for the Euler characteristic $\chi(\FF)$ of a rank $n$ coherent sheaf $\FF$ on the projective space $\mathbb P^N$ on the field $\mathbb K$, in terms of the Chern classes $c_i(\FF)$ of $\FF$.

For fixed $N$ (the dimension of the projective space) and $n$ (the rank of the sheaf), we explicitly obtain the polynomial $P(C_1, \dots, C_N)\in \mathbb  Q[C_1, \dots, C_N]$ such that
\[
P(c_1(\FF),\dots,c_N(\FF))=\chi(\FF)=\sum_{j=0}^N(-1)^j h^j\FF
\]
where we consider the Chern classes $c_i(\FF)$ on $\mathbb P^N$ as integers and  $h^j\FF$ is the dimension of the $i$-th cohomology module $H^i(\FF)$ as a $\mathbb K$-vector space.

More generally we obtain a general polynomial expression for  the Euler cha\-rac\-te\-ri\-stic of every twist of $\FF$ in terms of the Chern class $c_i(\FF)$ of $\FF$ and on $t$, namely:
\[
G(c_1(\FF),\dots,c_N(\FF),t)=\chi(\FF(t))=\sum_{j=0}^N(-1)^j h^j\FF(t)
\]
It is well known that such  polynomials exist and that they do not depend on $\FF$ (see for instance \cite{stable}, Theorem 2.3), so that they can be computed by means of the special cases given by  the totally split sheaves 
$\oplus \mathcal O_{\mathbb P^N}(a_i)$ with $a_i\geq 0$:  this argument is usually called the \lq\lq splitting principle\rq \rq. 

For our purposes those free sheaves  are very easy to manage  because their Euler cha\-rac\-te\-ri\-stic is only given by the 0-cohomology, $\chi(\oplus\OPN(a_i))=h^0(\oplus\OPN(a_i))$,  so that  it can be easily expressed as a sum of $n$ binomials involving $N$ and the $a_i$'s;  moreover also their Chern classes can be easily written  in terms of the $a_i$'s,  as   symmetric functions $c_j(\oplus\OPN(a_i))=\sum a_{i_1} \dots  a_{i_j}$, where the sum is over every sequence of $j$ indexes $i_1<\dots <i_j$.

However in practice a general computation necessary  involves \lq\lq changes of basis\rq\rq\ for polynomials,  mainly for those which are invariant  with respect to the action of  the permutation group on the variables,  that are not so easy to express in a suitable way. In fact we have to expand   polynomials expressed through binomials into their expansion as a sum of monomials and then  as a sum of elementary symmetric functions.
\bigskip

We will divide the solution of the problem in three steps.
\begin{enumerate}
	\item  We first write the polynomial  $Q_N(x_1, \dots, x_n )\in\mathbb Q[x_1,\dots,x_n]$ such that  $\chi(\oplus\OPN(a_i))= h^0\left(\oplus\OPN(a_i)\right)=\frac{1}{N!}Q_N(a_1, \dots,a_n)$ for every $a_i\geq 0$. As $a_1, \dots, a_n$ is not an ordered sequence, the polynomial $Q_N$ must be invariant under the action of the permutation group on the variable $x_i$.
	\item  Then we  \lq\lq change the basis\rq\rq\,  that is  we substitute the variables $x_i$ (corresponding to the $a_i$'s) by the variables $C_j$ given by  their  symmetric functions (corresponding to the Chern classes $c_i$) so obtaining the polynomial $P(C_1,\dots,C_N)$.
	
\end{enumerate}

The first step involves the Stirling numbers of first kind (that we introduce in \S\ref{stirling}); here we show that they give the coefficients of the expansion of a the polynomial $R_N(x)$ such that $h^0\OPN(a)=1/(N!)R_N(a)$ (see Theorem \ref{combin} and Corollary \ref{corollario}).

The second step is closely related to the {\it umbral calculus} (see \cite{rota} for an overview of the subject).
We also use the well-known {\it Newton-Girard formulas} (see \cite{seroul}) in order to obtain a faster algorithm to compute the polynomial $P$.

Finally using the relations between the Chern classes of a sheaf $\FF$ and those of the twists $\FF(t)$ (see (\ref{cherncontwist})), we compute the polynomial $G(C_1,\dots, C_N, T)\in \mathbb Q[C_1,\dots,C_N,T]$ such that $\chi(\FF(t))=G(c_1(\FF), \dots, c_N(\FF),t)$.

Beyond the theoretical results, we also present some procedures for the explicit computations of the polynomials $P(C_1,\dots,C_N)$ and $G(C_1,\dots,C_N,t)$, for a fixed dimension of the projective space $N$ and a fixed rank $n$ for the sheaf (see \S\ref{algoritmi}).

\section{Stirling Numbers of first kind}\label{stirling}

The {\it Stirling number of first kind} $\stir{N}{m}$ is the number of permutations of $N$ elements which contain exactly $m$ distinct cycles. \\
As a direct consequence of the definition, one can immediately see that: 
\begin{itemize}
	\item $\displaystyle \stir{0}{0}=1$ but $\displaystyle \stir{0}{m}=0 $ if  $m>0 $ and $\displaystyle  \stir{N}{0}=0 $ if $N>0$.
	
		\item $\displaystyle \stir{N}{m}=0 $  if $m>n $  and  $\displaystyle \stir{N}{N}=1$
		
	\item $\displaystyle \stir{N}{N-1}=\binom{N}{2}$ because a permutation of $N$ elements which contain $N-1$ cycles is determined by  its only 2-cycle.
\end{itemize}

Consider a \lq\lq square\rq\rq table whose entries are the integers $\stir{N}{m}$, where  each row is associated to a value for $N$ and each column is associated to a value for $m$. The above properties of  Stirling numbers  say that in such a table:
\begin{itemize}
	\item the triangle above the main diagonal is completely 0;
	\item on the main diagonal the entries are all 1's;
	\item on the \lq\lq second\rq\rq diagonal, the entries are the binomials $\binom{N}{2}$.
\end{itemize}

The following {\it recurrence relation} for the Stirling numbers of the first kind allows to complete the table:\\
\begin{equation}\label{ricorrenza}
\stir{N}{m}=\stir{N-1}{m-1}+(N-1)\stir{N-1}{m}
\end{equation}
This relation easily follow from the definition. In fact if we fix an element $\alpha$ among the $N$, there are two kinds of permutations of $N$ elements with $m$ cycles: the first addendum in 	the right side of (\ref{ricorrenza}) is the number of 
permutations containing  $(\alpha)$ as a 1-cycle; the second one 
corresponds to permutations not containing the 1-cycle $(\alpha) $: for every permutation of the other $N-1$ elements with $m$ cycles, the element $\alpha$   can be introduced in  $(N-1)$ different ways (in fact there are $j$ different ways to put a new element in a cycle of $j$ elements).

Now we can complete the table of the Stirling numbers of first kind
\begin{center}
		\begin{tabular}{|c|c c c c c c c c c c|}
			\hline
			 	& 0 & 1 & 2 & 3 & 4 & 5 & 6 & 7 & 8 & \dots\\
			 	\hline
			0 & 1 &		&		&		&		&		&		&  	& &\\
			1 & 0	& 1 & 	& 	& 	& 	& 	& 	& &\\
			2 &	0 &	1	&	1	&		&		&		&		&		& &\\
			3 &	0	&	2	&	3	& 1	&		&		&		&		& &\\
			4 &	0	&	6	&	11& 6	& 1	&		&		&		& &\\
			5 &	0	&	24&	50& 35& 11& 1	& 	&		& &\\
			6 &	0	&120&274&225& 85& 15& 1	&		& &\\
			7 &	0	&720&1764&1624&735&175&21& 1	&  &\\
			8 &	0	&5040&13068&13132&6769&1960&322&28& 1&\\
			\vdots &\vdots & \vdots& \vdots& \vdots& \vdots& \vdots& \vdots& \vdots&\vdots & $\ddots$\\
			\hline
		\end{tabular}
\end{center}

The above considered Stirling numbers of the first kind  are often called {\it unsigned}, in opposition to the {\it signed} Stirling numbers of first kind, which we denote by $s(N,m)$ and that are simply recovered from the unsigned ones by the rule
\[
s(N,m)=(-1)^{N-m}\stir{N}{m}.
\]
The original definition of signed Stirling numbers of first kind comes from a particular polynomial, the {\it falling factorial} $(x)_N$, which is ``similar'' to the one we are interested in:
\[
(x)_N=x(x-1)\cdots(x-N+1)
\]
The signed Stirling numbers of first kind are defined as the coefficients of the expansion
\[
(x)_N:=\sum_{k=0}^N s(N,k)x^k.
\]
We are mainly interested in the  {\it  rising factorial polynomial}
\[
(x)^{(N)}:=x(x+1)\cdots(x+N-1).
\]
or, more precisely to the polynomial \[
R_N(x):=(x+1)(x+2)\cdots(x+N) =(x+1)^{(N)}.
\]

\begin{Theorem}\label{combin} 
\begin{equation}\label{formula1}
R_N(x)=(x+1)^{(N)}=\sum_{k=0}^N\stir{N+1}{k+1}x^k.
\end{equation}
\end{Theorem}

\begin{pf}
We proceed by induction on $N$.

If $N=1$, the thesis is immediately verified.

Then assume that the formula holds for $R_{N-1}(x)$. Applying the  inductive hypothesis
to
\[
R_N(x)=(x+N)R_{N-1}(x)=(x+N)\left((x+1)^{(N-1)}\right)
\]
we obtain:
\[
R_N(x)=(x+N)\left(\sum_{k=0}^{N-1}\stir{N}{k+1}x^k\right)=\sum_{k=0}^{N-1}\stir{N}{k+1}x^{k+1}+\sum_{k=0}^{N-1}N\stir{N}{k+1}x^k.
\]
The change  $k+1 \rightarrow k$ in the first sum of the right side and the recurrence relation (\ref{ricorrenza}) give:
\begin{multline*}
R_N(x)=\stir{N}{N}x^N +\sum_{k=1}^{N-1}\left(\stir{N}{k}+N\stir{N}{k+1}\right)x^k+N\stir{N}{1}=\\
=\stir{N+1}{N+1}x^N +\sum_{k=1}^{N-1}\stir{N+1}{k+1}x^k+\stir{N+1}{1}
\end{multline*}
which is (\ref{formula1}) (note  that for every $r>0$,  $\stir{r}{r}=1$).  
\end{pf}

\begin{remark}
For a different  proof of Theorem \ref{combin}, we could refer to \cite{uff}, formula (6.11), and use the equality $x^{(N+1)}=x R_N(x)$.
\end{remark}

\begin{Corollary}\label{corollario}
For every $a\geq 0$ the dimension of the vector space of the degree $a$ hypersurfaces in $\mathbb{P}^N$ is given by:
\[
h^0\OPN(a)=\frac{1}{N!}R_N(a)=\frac{1}{N!}\sum_{k=0}^N\stir{N+1}{k+1}a^k.
\]
\end{Corollary}


\section{Invariant Polynomials and the Main Theorem}
Now let  $x_1, \dots, x_n$ be $n$ variables and consider the polynomial $Q_N(x_1, \dots,x_n):= \sum_j R_N(x_j) \in \mathbb{Q}[x_1, \dots, x_n]$.   This polynomial is closely  related to our problem, because  for every   choice of $n$  positive integers $a_1, \dots, a_n $, we have $\chi(\oplus \OPN(a_i))=h^0(\oplus \OPN(a_i))=\frac{1}{N!}Q_N(a_1, \dots,a_n)$. 

It is quite evident that $Q_N$ does not change under permutation of the variables, that is it is  \textit{invariant} for symmetric group $\mathcal S_n$. 
We just recall some basic definition and properties of the invariant polynomials; for more details one can see, for instance, \cite{goodman}.

The action of the symmetric group $\mathcal S_n$ on $\mathbb Q[x_1,\dots,x_n]$ is given in the following way.\\
If $p(x_1,\dots,x_n)\in \mathbb Q[x_1,\dots,x_n]$ and  $\sigma \in \mathcal S_n$, then:
\[
(\sigma \cdot p)(x_1,\dots,x_n)=p(x_{\sigma(1)},\dots,x_{\sigma(n)}).
\]

We say that $p\in \mathbb Q[x_1,\dots,x_n]$ is {\it invariant} for the action of $\mathcal S_n$ if
\[
\sigma \cdot p=p \quad \forall \sigma \in \mathcal S_n.
\]
It is easy to prove that the set of invariant polynomials that we denote by $\mathbb Q[x_1,\dots,x_n]^{\mathcal S_n}$ is an algebra, called  the {\it algebra of symmetric polynomials}.

Since $\mathcal S_n$ is a reductive linear algebraic group (see \cite{goodman}), there is a set of algebraically independent polynomials $\{f_1,\dots,f_n\}$, $f_i\in \mathbb Q[x_1,\dots,x_n]^{\mathcal S_n}$, such that  the polynomial ring they generate on $\mathbb Q$ is exactly   $\mathbb Q[x_1,\dots,x_n]^{\mathcal S_n}$, that is
\[
\mathbb Q[x_1,\dots,x_n]^{\mathcal S_n}=\mathbb Q[f_1,\dots,f_n].
\]
We call $\{f_1,\dots,f_n\}$ a set of {\it basic invariants}.

There are of course many sets of basic invariants for $\mathbb Q[x_1,\dots,x_n]^{\mathcal S_n}$, but we will be interested only in two of these:
\begin{itemize}
	\item the elementary symmetric polynomials: 
	\[
	C_0=1; \quad C_j:=\sum_{\lambda_1<\cdots<\lambda_j}x_{\lambda_1}\cdots  x_{\lambda_j} \quad j=1,\dots,n.
	\]
	\item the power sum symmetric polynomials:
	\[
	B_k:=\sum_{i=1}^n a_i^k \quad k=0,\dots, n.
	\]
\end{itemize}

Since both $\{C_1,\dots,C_n\}$ and $\{B_1,\dots,B_n\}$ are sets of basic invariants and so their elements  are algebraically independent, we can consider them as indeterminates.

 Every invariant polynomial, included $Q_N(x_1, \dots,x_n)= \sum_i R_N(x_i)$, can be written using either of the two sets of basic invariants: as a polynomial in the indeterminates $C_j$'s and as a polynomial in the indeterminates $B_k$'s. 
\begin{Lemma}\label{primachi}
In the above notation:
\[
Q_N(x_1, \dots,x_n) =\sum_{k=0}^N\stir{N+1}{k+1}B_k.
\]
\end{Lemma}

\begin{pf}
Applying Theorem \ref{combin}, we immediately obtain
\[
Q_N(x_1, \dots,x_n)=\sum_{i=1}^n\left( \sum_{k=0}^N\stir{N+1}{k+1}x_i^k\right)=\sum_{k=0}^N \stir{N+1}{k+1}\left( \sum_{i=1}^nx_i^k\right).
\]
where we can commute the two summations because they are independent.
\end{pf}

If we know the expression of an invariant polynomial  in terms of a set of basic invariants and  want to obtain its expression in terms of the other one, we have to manage the not so easy problem of the \lq\lq change of basis\rq\rq. For instance:
\[
B_0=n C_0 \quad , \quad B_1=C_1 \quad , \quad B_2=C_1^2-2C_2.
\]
In order to find a general  expression of $B_k$ as a function of the $C_j$'s, we recall the Newton-Girard formula (see \cite{seroul}):
\[
(-1)^rB_r+\sum_{l=1}^r(-1)^{r+l}B_l C_{r-l}=0
\]
Note that in fact this formula    holds for every $r\in \mathbb{N}$, with the convention that $B_k$  is the sum of powers $x_i^k$ and $C_k=0$ if $k\geq n+1$.

With these notations, we can then prove
\begin{Lemma}\label{cambio}
For every $1\leq r\leq n$ $B_r=\det (M_r)$ where
\begin{equation} \label{matrice}
				M_r=\left(
				\begin{array}{ccccccc}
				C_1 	&1 		& 0 	& 0 & \ldots & 0 & 0	 \\
				2C_2 	&C_1	& 1		& 0	& \ldots & 0 & 0	 \\
				3C_3	&C_2	&C_1	& 1	& \ldots & 0 & 0  \\
				\vdots & \vdots & \vdots &\vdots & & \vdots & \vdots \\
				(r-1)C_{r-1} & C_{r-2} & C_{r-3} & C_{r-4}&\ldots & C_1 & 1\\
				rC_r & C_{r-1} & C_{r-2} & C_{r-3}&\ldots & C_2 & C_1
				\end{array} \right)
\end{equation}
\end{Lemma}

\begin{pf}
We proceed by induction on $r$.
If $r=1$ there is nothing to prove.

Assume $r\geq 2$ and the thesis true for $B_l$, $l\leq r-1$.
From Newton-Girard formula we have
\begin{equation}\label{det}
B_r=B_{r-1}C_1-B_{r-2}C_2+\dots+(-1)^{r-2}B_1C_{r-1}+(-1)^{r-1}rC_r.
\end{equation}
Observe that since the thesis is true for $B_l$, $l\leq r-1$, we can write the second term of ($\ref{det}$) as
\begin{equation}\label{ricorsiva}
B_r=\sum_{l=1}^{r-1} (-1)^{l-1}C_l\det (M_{r-l})+(-1)^{r-1}rC_r.
\end{equation}
Thanks to the presence of the $1$'s above the main diagonal of the matrix (\ref{matrice}), this is exactly the determinant of $M_r$.
\end{pf}

\begin{remark}
For a different proof of the previous result one can see  \cite{little}, Chapter $V$. However the proof we present is more constructive and gives rise to a faster algorithm, that we will present in \S \ref{algoritmi}.
\end{remark}

Finally, we obtain the main Theorem as an application of Lemmas \ref{primachi} and \ref{cambio}

\begin{Theorem}\label{finalmente}
Let $\FF$ be rank $n$  reflexive sheaf on $\mathbb P^N$. Consider
\[
P=\frac{1}{N!}\sum_{k=1}^N\stir{N+1}{k+1}\det (M_k)+n
\]
where  $M_k$ is the $k \times k$ matrix previously defined and its determinant is a polynomial in the variables $C_1,\dots,C_k$.\\
Then
\[P(c_1(\FF),\dots,c_N(\FF))=\chi(\FF)\].
\end{Theorem}

\begin{pf}
First, using the \lq\lq splitting principle\rq\rq, we know there is a polynomial  $ P \in \mathbb Q[c_1,\dots,c_N]$, depending only on $N$ and $n$, such that
$\chi(\FF)=P(c_1(\FF), \dots,c_N(\FF))$ for any rank $n$ coherent sheaf $\FF$.

It is then sufficient to find such a polynomial for the sheaf $\oplus_{i=1}^n\OPN(a_i)$, $a_i\geq 0$, with $n\geq N$.

Thanks to Lemma \ref{cambio}, we can pass from the $B_k$'s to the $C_j$'s in the expression of Lemma \ref{primachi} for $Q_N(x_1,\dots,x_n)$:
\[
Q_N(x_1,\dots,x_n)=\sum_{k=0}^N\stir{N+1}{k+1}B_k=\sum_{k=1}^N\stir{N+1}{k+1}\det (M_k)+n(N!).
\]

Using Corollary \ref{corollario}, we obtain
\begin{multline*}
\chi(\oplus_{i=1}^n\OPN(a_i))=h^0(\oplus_{i=1}^n\OPN(a_i))=\frac{1}{N!}Q_N(a_1,\dots,a_n)=\\
=\frac{1}{N!}\sum_{k=1}^N\stir{N+1}{k+1}\det (M_k)(c_1,\dots,c_k)+n
\end{multline*}
where $\det (M_k)(c_1,\dots,c_k)$ means evaluating the polynomial in the Chern classes of $\oplus_{i=1}^n\OPN(a_i)$.

\end{pf}

\begin{remark}
In the proof of Theorem \ref{finalmente}, the assumption $n\geq N$ for the completely split bundle is not a lost in generality; in fact, a coherent sheaf $\FF$ of rank $n\leq N$ may have  $c_i(\FF)\neq 0$ for $i\geq n+1$.
\end{remark}

\section{Implementation}\label{algoritmi}

In the previous paragraphs we obtained the  following result:

\begin{quote}
	Let $\FF$ be a rank $n$ coherent sheaf on $\mathbb P^N$. Then 
	\[
	\chi(\FF)=\frac{1}{N!}\sum_{k=1}^N\stir{N+1}{k+1}\det M_k(\FF)+n
	\]
	where $M_k$ is a $k \times k$ matrix whose definition is (\ref{matrice}), $\det(M_k)\in \mathbb Z[c_1,\dots,c_k]$ and with $\det M_k(\FF)$ we mean $\det M_k(c_1(\FF),\dots,c_N(\FF)]$.
\end{quote}

The polynomial for $\chi(\FF)$ is not too easy to handle, since it contains some determinants.

Anyway, if we fix the dimension $N$, it is quite easy to write a procedure to compute the polynomial $\chi(\FF)$.

\medskip

Here we write some procedures for Maple. Probably they are not the best implementations for the algorithms we wish to expose, they are just intended to be examples.

\medskip

\subsection{A first algorithm for $\chi(\FF)$}
First we write a procedure to write the $r$-th row of the matrix $M_k$

\begin{verbatim}
 Row:=proc(r,k) 
 v:=[c[1]];
 if (r=1) 
   then v:=[op(v),1];
       for j from 2 to n-r do v:=[op(v),0] od; 
       return v; 
    else 
     if (r=2) then 
 		     if (k=2) 
           then return  [2*c[2], c[1]]; 
           else v:=[2*c[2],c[1],1]  fi; 
     for j from 2 to k-r do v:=[op(v),0] od; 
     return v ;
     else 
       for i from 2 to r-1 do v:=[c[i],op(v)] od 
     fi 
 fi ;
 v:=[r*c[r],op(v)]; 
 if (r=k) 
 	  then return v; 
 	  else v:=[op(v),1]
 fi;
 for j from 2 to k-r do v:=[op(v),0] od; 
 return v;
 end proc;
\end{verbatim}

Then we write the procedure that outputs the matrix $M_k$

\begin{verbatim}
 MatrixM:=proc(k)
 if (k=1) 
    then return matrix(1,1,[c[1]]) 
 fi; 
 V:=[];
 for i from 1 to k do V:=[op(V),op(Row(i,k))] 
 od;
 return matrix(k,k,[op(V)]);
 end proc;
\end{verbatim}

Finally, we write the procedure that returns the polynomial for $\chi(\FF)$ once that we have fixed $N$

\begin{verbatim}
 chi:=proc(n,N) 
 for i from 1 to N do S[i]:=linalg[det](MatrixM(i)) 
 od;
 return 1/N!*sum(abs(stirling1(N+1,k+1))*S[k],k=1..N)+n; 
 end proc;
\end{verbatim}

With the procedure \verb+chi(n,N)+ one can then easily obtain the Euler characteristic for a rank $n$ coherent sheaf on $\mathbb P^N$, just fixing $N$.

The polynomial expression for the  Euler Characteristic for a rank $n$ sheaf on $\mathbb P^3$ is known: one can see \cite{stable}, Theorem 2.3. We give as examples the polynomial expressions for $\chi(\FF)$ for a rank $n$ sheaf on $\mathbb P^4$ and $\mathbb P^5$.

\medskip

$chi(n,4)$;
\[
\frac{1}{24}\left[{c_{{1}}}^{4}+10\,{c_{{1}}}^{3}-4\,{c_{{1}}}^{2}c_{{2}}+35\,{c_{{1}}}^
{2}-30\,c_{{1}}c_{{2}}+4\,c_{{1}}c_{{3}}+\right.\]
\[\left.2\,{c_{{2}}}^{2}+50\,c_{{1}}-
70\,c_{{2}}+30\,c_{{3}}-4\,c_{{4}}\right]+n\]
 
$chi(n,5)$;
\[
\frac{1}{120}\left[{c_{{1}}}^{5}+15\,{c_{{1}}}^{4}-5\,{c_{{1}}}^{3}c_{{2}}+85\,{c_{{1}}}^
{3}-60\,{c_{{1}}}^{2}c_{{2}}+5\,{c_{{1}}}^{2}c_{{3}}+5\,c_{{1}}{c_{{2}
}}^{2}+225\,{c_{{1}}}^{2}+\right.\]
\[-255\,c_{{1}}c_{{2}}+60\,c_{{1}}c_{{3}}-5\,c_
{{1}}c_{{4}}+30\,{c_{{2}}}^{2}-5\,c_{{2}}c_{{3}}+274\,c_{{1}}+\]
\[\left.-450\,c_{
{2}}+255\,c_{{3}}-60\,c_{{4}}+5\,c_{{5}}\right]+n
\]

\subsection{A faster algorithm for $\chi(\FF)$}
The procedure \verb+chi(n,N)+ is quite expensive  from the computational viewpoint: Maple 11 on a personal computer (Intel Pentium CPU 3.00 Ghz, 992 mb RAM) took more than 20 seconds for the case $N=18$. 

We can improve the procedure because actually we do not need to construct the matrices $M_k$ to compute their determinant. We can just construct a recursive procedure using formula (\ref{ricorsiva}).

First the procedure to compute $\det M_k$:
\begin{verbatim}
detM:=proc(k) 
if k=1 then return c[1] 
fi;
if k=2 then return c[1]^2-2*c[2] 
fi; 
M:=(-1)^(k-1)*k*c[k]; 
for i from 1 to k-1 do M:=M+(-1)^(i-1)*c[i]*detM(k-i) 
od;
return expand(M) end proc;
\end{verbatim}

Then we rewrite the procedure \verb+chi(n,N)+, but using \verb+detM(k)+:

\begin{verbatim}
 chifast:=proc(n,N) 
 for i from 1 to N do M[i]:=detM(i) 
 od;
 return 1/N!*sum(abs(stirling1(N+1,k+1))*M[k],k=1..N)+n; 
 end proc;
\end{verbatim}

This last procedure is much faster than \verb+chi(n,N)+: for instance, we computated the polynomial for $\chi(\FF)$ for a rank $n$ coherent sheaf $\FF$ on $\mathbb P^{20}$, on a personal computer (Intel Pentium CPU 3.00 Ghz, 992 mb RAM) using Maple 11:
\begin{itemize} 
	\item \verb+chi(n,20)+ took 172.14 sec to output the polynomial;
	\item \verb+chifast(n,20)+ took only 7.72 sec to output the polynomial.
\end{itemize}

\subsection{An algorithm for $\chi(\FF(t))$}

Since we have already a polynomial form for $\chi(\FF)$, we can easily obtain the polynomial associated to $\chi(\FF(t))$ for every $t\in\mathbb Z$.
It is sufficient to remember that, if $\FF$ is a rank $n$ coherent sheaf on $\mathbb P^N$ and Chern classes $c_i$, then
\begin{equation}\label{cherncontwist} 
c_i(\FF(t))=c_i+(n-i+1)tc_{i-1}+{n-i+2\choose 2}t^2c_{i-2}+\cdots+{n-1\choose i-1}t^{i-1}c_1+{n\choose i}t^i. 
\end{equation}

So we substitue $C_i(T)=C_i+(n-i+1)TC_{i-1}+{n-i+2\choose 2}T^2C_{i-2}+\cdots+{n-1\choose i-1}T^{i-1}C_1+{n\choose i}T^i$ to $C_i$ in $P(C_1,\dots,C_N)$  obtaining
\[
P(C_1(T),\dots,C_N(T))=G(C_1,\dots,C_N,T)\in \mathbb Q[C_1,\dots,C_N,T].
\]

With some little changes, we can rewrite the procedure \verb+chifast(n,N)+ for any twist $\FF(t)$, $t\in \mathbb Z$: the procedure outputs a polynomial in the variable $T$.
First, we write a procedure to obtain a \lq\lq twisted\rq\rq Chern class
\begin{verbatim}
ct:=proc(j,N) 
cT:=c[j];
for i from 1 to j-1 do
cT:=cT+binomial(N-j+i,i)*T^i*c[j-i] 
od;
cT:=cT+binomial(N,j)*T^j;
return cT;
end proc;
\end{verbatim}

Then we simply rewrite the procedure \verb+detM(k)+ 
\begin{verbatim}
detMT:=proc(k,N)
if (k=1) then return ct(1,N) 
fi;
if k=2 then return ct(1,N)^2-2*ct(2,N) 
fi;
Mt:=(-1)^(k-1)*k*ct(k,N); 
for i from 1 to k-1 do Mt:=Mt+(-1)^(i-1)*ct(i,N)*detMt(k-i,N) 
od;
return expand(Mt) 
end proc;
\end{verbatim}

Finally, we rewrite \verb+chifast(r,N)+ using the \lq\lq twisted\rq\rq determinants
\begin{verbatim}
chit:=proc(n,N)
for i from 1 to N do M[i]:=detMT(i,N)
od;
return sort(collect(1/N!*sum(abs(stirling1(N+1,k+1))*M[k],k=1..N)+n,T),T); 
end proc;
\end{verbatim}

For instance, we obtain the polynomial associated to $\chi(\FF(t))$ for a coherent sheaf on $\mathbb P^6$\\

$chit(n,6)$;
\[
\frac{1}{240}\left[6\,{T}^{6}+ \left( 6\,c_{{1}}+126 \right) {T}^{5}+ \left( 15\,{c_{{1}}}^{2}+105\,c_{{1}}+1050-30\,c_{{2}} \right) {T}^{4}+\right.\]
\[+ \left( 60\,c_{{3}}+210\,{c_{{1}}}^{2}+20\,{c_{{1}}}^{3}+700\,c_{{1}}+4410-60\,c_{{1}}c_{{2}}-420\,c_{{2}} \right){T}^{3}+\]
\[+\left(-60\,c_{{4}}+1050\,{c_{{1}}}^{2}-60\,{c_{{1}}}^{2}c_{{2}}-2100\,c_{{2}}-630\,c_{{1}}c_{{2}}+2205\,c_{{1}}+\right.\]
\[\left.+630\,c_{{3}}+9744+210\,{c_{{1}}}^{3}+60\,c_{{1}}c_{{3}}+15\,{c_{{1}}}^{4}+30\,{c_{{2}}}^{2}\right){T}^{2}+\]
\[+\left(420\,c_{{1}}c_{{3}}-4410\,c_{{2}}-30\,{c_{{1}}}^{3}c_{{2}}+30\,c_{{1}}{c_{{2}}}^{2}+2205\,{c_{{1}}}^{2}+2100\,c_{{3}}+30\,c_{{5}}+\right.\]
\[+700\,{c_{{1}}}^{3}+105\,{c_{{1}}}^{4}+6\,{c_{{1}}}^{5}+3248\,c_{{1}}+210\,{c_{{2}}}^{2}-420\,c_{{4}}+10584-30\,c_{{2}}c_{{3}}+30\,{c_{{1}}}^{2}c_{{3}}+\]
\[\left.-30\,c_{{1}}c_{{4}}-2100\,c_{{1}}c_{{2}}-420\,{c_{{1}}}^{2}c_{{2}}\right)T+2205\,c_{{3}}-3248\,c_{{2}}-700\,c_{{4}}+r+350\,{c_{{2}}}^{2}+\]
\[+175\,{c_{{1}}}^{4}+1764\,c_{{1}}-105\,c_{{1}}c_{{4}}+105\,{c_{{1}}}^{2}c_{{3}}-105\,{c_{{1}}}^{3}c_{{2}}+105\,c_{{5}}-105\,c_{{2}}c_{{3}}+1624\,{c_{{1}}}^{2}+\]
\[-2205\,c_{{1}}c_{{2}}+105\,c_{{1}}{c_{{2}}}^{2}+700\,c_{{1}}c_{{3}}-700\,{c_{{1}}}^{2}c_{{2}}-6\,{c_{{1}}}^{2}c_{{4}}+6\,{c_{{1}}}^{3}c_{{3}}-6\,{c_{{1}}}^{4}c_{{2}}+9\,{c_{{1}}}^{2}{c_{{2}}}^{2}+\]
\[\left.+{c_{{1}}}^{6}-2\,{c_{{2}}}^{3}+6\,c_{{2}}c_{{4}}+3\,{c_{{3}}}^{2}+21\,{c_{{1}}}^{5}-6\,c_{{6}}+735\,{c_{{1}}}^{3}-12\,c_{{1}}c_{{2}}c_{{3}}+6\,c_{{1}}c_{{5}}\right]+n.
\]

\end{document}